%% file: CRNormal_no.tex
\documentclass[10pt,reqno]{amsart}
\usepackage{amssymb,amscd,amsbsy}
\usepackage{amsthm}
\input{CRinput}

\theoremstyle{remark}

\newtheorem*{rem*}{Remark}

\begin{document}

\newcommand{\vse}{\vspace{.2in}}

\title{Functions  of perturbed normal operators}

\maketitle
\begin{center}
\Large
Aleksei Aleksandrov$^{\rm a}$, Vladimir Peller$^{\rm b}$, Denis Potapov$^{\rm c}$, Fedor Sukochev$^{\rm c}$
\end{center}

\begin{center}
\footnotesize
{\it$^{\rm a}$St-Petersburg Branch, Steklov Institute of Mathematics, Fontanka 27, 191023 St-Petersburg, Russia\\
$^{\rm b}$Department of Mathematics, Michigan State University, East Lansing, MI 48824, USA\\
$^{\rm c}$School of Mathematics \& Statistics, University of NSW, Kensington NSW 2052, Australia}
\end{center}

\newcommand{\mt}{{\mathcal T}}

\footnotesize

{\bf Abstract.}
In \cite{Pe1}, \cite{Pe2}, \cite{AP1}, \cite{AP2}, and \cite{AP3} sharp estimates for $f(A)-f(B)$ were obtained for self-adjoint operators $A$ and $B$ and for various classes of functions $f$ on the real line $\R$.
In this note we extend those results to the case of functions of normal operators. We show that if $f$ belongs to the H\"older class $\L_\a(\R^2)$, $0<\a<1$, of functions of two variables, and $N_1$ and $N_2$ are normal operators, then
$\|f(N_1)-f(N_2)\|\le\const\|f\|_{\L_\a}\|N_1-N_2\|^\a$. We obtain a more general result for functions in the space
$\L_\o(\R^2)=\big\{f:~|f(\z_1)-f(\z_2)|\le\const\o(|\z_1-\z_2|)\big\}$ for an arbitrary modulus of continuity $\o$. We prove that if $f$ belongs to the Besov class $B_{\be1}^1(\R^2)$, then it is operator Lipschitz, i.e.,
$\|f(N_1)-f(N_2)\|\le\const\|f\|_{B_{\be1}^1}\|N_1-N_2\|$. We also study properties of $f(N_1)-f(N_2)$ in the case when $f\in\L_\a(\R^2)$ and $N_1-N_2$ belongs to the Schatten-von Neuman class $\bS_p$.

\medskip

\begin{center}
{\bf\large Fonctions d'op\'erateurs perturb\'es normaux}
\end{center}

\medskip

{\bf R\'esum\'e.} On a obtenu dans \cite{Pe1}, \cite{Pe2}, \cite{AP1}, \cite{AP2} et \cite{AP3} des estimations pr\'ecises de $f(A)-f(B)$, o\`u $A$ et $B$ sont des op\'erateurs autoadjoints et $f$ est une fonction sur la droite r\'eelle $\R$.  Dans cette note nous obtenons des g\'en\'eralisations de ces r\'esultats pour les op\'erateurs normaux et pour les fonctions $f$ de deux variables. Nous
d\'emontrons que si $f$ appartient \`a l'espace de H\"older $\L_\a(\R^2)$, $0<\a<1$, alors
$\|f(N_1)-f(N_2)\|\le\const\|f\|_{\L_\a}\|N_1-N_2\|^\a$ pour chaque op\'erateurs normaux $N_1$ et $N_2$. Nous obtenons aussi un r\'esultat plus g\'en\'eral pour les fonctions de la classe
$\L_\o(\R^2)=\lb\big\{f:~|f(\z_1)-f(\z_2)|\le\const\o(|\z_1-\z_2|)\big\}$. Nous montrons que si $f$ appartient \`a l'espace de Besov $B_{\be1}^1(\R^2)$, alors
$f$ est une fonction lipschitzienne op\'eratorielle, c'est-\`a-dire
$\|f(N_1)-f(N_2)\|\le\const\|f\|_{B_{\be1}^1}\|N_1-N_2\|$ pour chaque op\'erateurs normaux $N_1$ et $N_2$. Nous aussi \'etudions les propri\'et\'es de
$f(N_1)-f(N_2)$ en cas quand $f\in\L_\a(\R^2)$ et $N_1$ et $N_2$ sont des op\'erateurs normaux tells que $N_1-N_2$ appartient \`a l'espace $\bS_p$ de Schatten--von Neumann.

\normalsize

\

\begin{center}
{\bf\large Version fran\c caise abr\'eg\'ee}
\end{center}

\medskip

Il est bien connu (voir \cite{F}) qu'il y a des fonctions $f$ lipschitziennes sur la droite r\'eelle $\R$ qui ne sont pas {\it lipschitziennes op\'eratorielles},
c'est-\`a-dire la condition
$$
|f(x)-f(y)|\le\const|x-y|,\quad x,~y\in\R,
$$
n'implique pas que pour tous les op\'erateurs auto-adjoints $A$ et $B$ l'in\'egalit\'e
$$
\|f(A)-f(B)\|\le\const\|A-B\|.
$$
soit vraie.

Dans \cite{Pe1} et \cite{Pe2} des conditions n\'ecessaires et des conditions suffisantes \'etaient trov\'ees pour qu'une fonction $f$ soit lipschitzienne
op\'eratorielle. En particulier, il \'etait d\'emontr\'e dans \cite{Pe1}
que pour qu'une fonction $f$ soit lipschitzienne op\'eratorielle il est n\'ecessaire que $f$ appartienne locallement \`a l'espace de Besov $B_{11}^1(\R)$.
C'est aussi implique qu'une fonction lipschitzienne n'est pas n\'ecessairement
lipschitzienne op\'eratorielle.

D'autre part il \'etait d\'emontr\'e dans \cite{Pe1} et \cite{Pe2} que si $f$ appartient dans l'espace de Besov $B_{\be1}^1(\R)$, alors $f$ est lipschitzienne
op\'eratorielle.

Il se trouve que la situation change dramatiquement si l'on consid\`ere les fonctions de la classe $\L_\a(\R)$ de H\"older d'ordre $\a$, $0<\a<1$.
Il \'etait d\'emontr\'e dans \cite{AP1} et \cite{AP2}  que si
$f$ apartient \`a $\L_\a(\R)$, $0<\a<1$,
(c\'est-\`a-dire $|f(x)-f(y)|\le\const|x-y|^\a$), alors $f$ doit \^etre
{\it h\"olderienne op\'eratorielle d'ordre} $\a$, c\'est-\`a-dire
$$
\|f(A)-f(B)\|\le\const\|A-B\|^\a.
$$

Dans \cite{AP1} et \cite{AP2} un probl\`eme plus g\'en\'eral \'etait consid\'er\'e pour les fonctions dans l'espace
$$
\L_\o(\R)=\big\{f:~|f(\z_1)-f(\z_2)|\le\const\o(|\z_1-\z_2|)\big\},
$$
o\`u $\o$ est un module de continuit\'e arbitraire.

Finalement il \'etait d\'emontr\'e dans \cite{AP1} et \cite{AP3} que si $A$ et $B$ sont des op\'erateurs autoadjoints tels que $A-B$ appartient \`a la classe de Schatten--von Neumann $\bS_p$, $p>1$, et $f\in\L_\a(\R)$, $0<\a<1$, alors
$f(A)-f(B)\in \bS_{p/\a}$.

Dans cette note nous g\'en\'eralisons les r\'esultats ci-dessus aux cas des op\'erateurs normaux (pas n\'ecessairement born\'es).

Nos r\'esultats sont bas\'es sur l'in\'egalit\'e suivante:
\bay
\label{nB}
\|f(N_1)-f(N_2)\|\le\const\s\|f\|_{L^\be}\|N_1-N_2\|,
\ey
o\`u $f$ est une fonction born\'ee sur $\R^2$ dont la transform\'ee de Fourier
a un support dans le disque $\big\{(x,y)\in\R^2:~x^2+y^2\le\s^2\big\}$
et $N_1$ et $N_2$ sont des op\'erateurs normaux.

Pour \'etablir \rf{nB}, nous utilisons les int\'egrales doubles op\'eratorielles
et nous
d\'emontrons la formule suivante:
\begin{align*}
f(N_1)-f(N_2)&=\iint\limits_{\C^2}\frac{f(x_1,y_1)-f(x_1,y_2)}{y_1-y_2}\,
dE_1(z_1)(B_1-B_2)\,dE_2(z_2)\\[.2cm]
&+\iint\limits_{\C^2}\frac{f(x_1,y_2)-f(x_2,y_2)}{x_1-x_2}\,
dE_1(z_1)(A_1-A_2)\,dE_2(z_2),
\end{align*}
o\`u $N_1$ et $N_2$ sont des op\'erateurs normaux dont la diff\'erence est born\'ee,
$E_1$ et $E_2$ sont les mesures sp\'ectrales d'$N_1$ et d'$N_2$. Ici
$x_j=\re z_j$, $y_j=\im z_j$, $A_j=\re N_j$, $B_j=\im N_j$, $j=1,2$.

En utilisant \rf{nB} nous d\'emontrons que si $f$ appartient \`a la classe de Besov $B_{\be1}^1(\R^2)$, alors $f$ est une fonction lipschitzienne op\'eratorielle et
$$
\|f(N_1)-f(N_2)\|\le\const\|f\|_{B_{\be1}^1(\R^2)}\|N_1-N_2\|
$$
pour chaque op\'erateurs normaux $N_1$ et $N_2$.

D'autre part, si $f$ appartient \`a la classe de H\"older $\L_\a(\R^2)$, $0<\a<1$, alors
$$
\|f(N_1)-f(N_2)\|\le\const\|f\|_{\L_\a(\R^2)}\|N_1-N_2\|^\a
$$
pour chaque op\'erateurs normaux $N_1$ et $N_2$.

Supposons maintenant que $\o$ est un module de continuit\'e et
$$
\o_*(x)=x\int_x^\be\frac{\o(t)}{t^2}\,dt,\quad x>0,
$$
alors
$$
\|f(N_1)-f(N_2)\|\le\const\|f\|_{\L_\o(\R^2)}\,\o_*\big(\|N_1-N_2\|\big)
$$
pour chaque op\'erateurs normaux $N_1$ et $N_2$.

Finalement nous d\'emontrons que si $f\in\L_\a(\R^2)$, $0<\a<1$, et $N_1$ et $N_2$ sont des op\'erateurs normaux dont la difference appartient \`a l'espace
de Schatten--von Neumann $\bS_p$, $p>1$, alors $f(N_1)-f(N_2)\in\bS_{p/\a}$ et
$$
\|f(N_1)-f(N_2)\|_{\bS_{p/\a}}\le\const\|f\|_{\L_\a(\R^2)}\|N_1-N_2\|_{\bS_p}^\a.
$$

\begin{center}
------------------------------
\end{center}

\setcounter{section}{0}
\section{\bf Introduction}

\medskip

In this note we generalize results of the papers \cite{Pe1}, \cite{Pe2},
\cite{AP1}, \cite{AP2}, and \cite{AP3} to the case of normal operators.

A Lipschitz function $f$ on the real line $\R$ (i.e., a function satisfying
the inequality $|f(x)-f(y)|\le\const|x-y|$, $x,\,y\in\R$) does not have to be {\it operator Lipschitz}, i.e,
$$
\|f(A)-f(B)\|\le\const\|A-B\|
$$
for arbitrary self-adjoint operators $A$ and $B$ on Hilbert space.
The existence of such functions was proved in \cite{F}.
Later in \cite{Pe1} and \cite{Pe2} necessary conditions were found for a function $f$ to be operator Lipschitz. In particular, it was shown in \cite{Pe1} that if $f$ is operator Lipschitz, then $f$ belongs locally to the Besov space $B_{11}^1(\R)$. This also implies that Lipschitz  functions do not have to be operator Lipschitz. Note that in \cite{Pe1} and \cite{Pe2} stronger necessary conditions are also obtained.
Note also that the necessary conditions obtained in \cite{Pe1} and \cite{Pe2} are based on the trace class criterion for Hankel operators, see \cite{Pe3}, Ch. 6.

On the other hand, it was shown in \cite{Pe1} and \cite{Pe2} that if $f$ belongs to the Besov class $B_{\be1}^1(\R)$, then $f$ is operator Lipschitz. We refer the reader to \cite{Pee} for information on Besov spaces.

It was shown in \cite{AP1} and \cite{AP2}
that the situation dramatically changes if we consider H\"older classes $\L_\a(\R)$ with $0<\a<1$. In this case such functions are necessarily {\it operator H\"older of order $\a$}, i.e., the condition
$|f(x)-f(y)|\le\const|x-y|^\a$, $x,\,y\in\R$,
implies that for self-adjoint operators
$A$ and $B$ on Hilbert space,
$$
\|f(A)-f(B)\|\le\const\|A-B\|^\a.
$$
Note that another proof of this result was found in \cite{FN}.

This result was generalized in \cite{AP1} and \cite{AP2} to the case of functions
of class $\L_\o(\R)$ for arbitrary moduli of continuity $\o$. This class consists of functions $f$ on $\R$, for which
$|f(x)-f(y)|\le\const\o(|x-y|)$, $x,\,y\in\R$,

Finally, we mention here that in \cite{AP3} properties of operators $f(A)-f(B)$
were studied for functions $f$ in $\L_\a(\R)$ and self-adjoint operators $A$ and $B$ whose difference $A-B$ belongs to Schatten--von Neumann classes $\bS_p$.

In this paper we generalize the above results to the case of (not necessarily bounded) normal operators. Throughout the paper we identify the complex plane 
$\C$ with $\R^2$.

\medskip

\section{\bf Double operator integrals and the key inequality}

\medskip

Our results are based on the following inequality:

\begin{thm}
\label{fn}
Let $f$ be a bounded function of class $L^\be(\R^2)$ whose Fourier transform is supported on the disc $\{\z\in\C:~|\z|\le\s\}$. Then
$$
\|f(N_1)-f(N_2)\|\le\const\s\,\|N_1-N_2\|
$$
for arbitrary normal operators $N_1$ and $N_2$ with bounded difference.
\end{thm}

To prove Theorem \ref{fn}, we obtain a formula for $f(N_1)-f(N_2)$ in terms of double operator integrals. The theory of double operator integrals was developed in \cite{BS1}, \cite{BS2}, and \cite{BS3}. If $E_1$ and $E_2$ are spectral measures on $\X_1$ and $\X_2$ and $\Phi$ is a bounded measurable function
on $\X_1\times\X_2$, then the double operator integral
$$
\iint\limits_{\X_1\times\X_2}\Phi(s_1,s_2)\,dE_1(s_1)T\,dE_2(s_2)
$$
is well defined for all operators $T$ of Hilbert--Schmidt class $\bS_2$ and determines an operator of class $\bS_2$. For certain functions $\Phi$ the transformer $T\mapsto \iint\Phi\,dE_1T\,dE_2$ maps the trace class $\bS_1$ into itself. For such functions $\Phi$ one can define by duality double operator integrals for all bounded operators $T$. Such functions $\Phi$ are called {\it Schur multipliers}
(with respect to the spectral measures $E_1$ and $E_2$). We refer the reader to \cite{Pe1} for characterizations of Schur multipliers.

In the following theorem $E_1$ and $E_2$ are the spectral measures of normal operators $N_1$ and $N_2$. We use the notation $x_j=\re z_j$, $y_j=\im z_j$, $A_j=\re N_j$, $B_j=\im N_j$, $j=1,2$.

\begin{thm}
\label{doi}
Let $N_1$ and $N_2$ be normal operators such that $N_1-N_2$ is bounded. Suppose that $f$ is a function in $L^\be(\R^2)$ such that its Fourier transform $\F f$ has compact support. Then
the functions
$$
(z_1,z_2)\mapsto\frac{f(x_1,y_1)-f(x_1,y_2)}{y_1-y_2}\qquad\mbox{and}\qquad
(z_1,z_2)\mapsto\frac{f(x_1,y_2)-f(x_2,y_2)}{x_1-x_2}
$$
are Schur multipliers with respect to $E_1$ and $E_2$ and
\begin{align}
\label{if}
f(N_1)-f(N_2)&=\iint\limits_{\C^2}\frac{f(x_1,y_1)-f(x_1,y_2)}{y_1-y_2}\,
dE_1(z_1)(B_1-B_2)\,dE_2(z_2)\nonumber\\[.2cm]
&+\iint\limits_{\C^2}\frac{f(x_1,y_2)-f(x_2,y_2)}{x_1-x_2}\,
dE_1(z_1)(A_1-A_2)\,dE_2(z_2).
\end{align}
\end{thm}

\medskip

\section{\bf Operator Lipschitz functions of two variables}
\setcounter{equation}{0}

\medskip

A continuous function $f$ on $\R^2$ is called {\it operator Lipschitz} if
$$
\|f(N_1)-f(N_2)\|\le\const\|N_1-N_2\|
$$
for arbitrary normal operators $N_1$ and $N_2$ whose difference is a bounded operator.

\begin{thm}
\label{oL}
Let $f$ belong to the Besov space $B_{\be1}^1(\R^2)$
and let $N_1$ and $N_2$ be normal operators whose difference is a bounded operator. Then {\em\rf{if}} holds and
$$
\|f(N_1)-f(N_2)\|\le\const\|f\|_{B_{\be1}^1(\R^2)}\|N_1-N_2\|.
$$
\end{thm}

In other words, functions in $B_{\be1}^1(\R^2)$ must be operator Lipschitz.

As in the case of functions on $\R$, not all Lipschitz functions are operator Lipschitz. In particular, it follows from \cite{Pe1} that if $f$ is an operator Lipschitz function on $\R^2$, then the restriction of $f$ to an arbitrary line belongs locally to the Besov space $B_{11}^1$.

The next result shows that functions in $B_{\be1}^1(\R^2)$ respect trace class perturbations.

\begin{thm}
\label{yad}
Let $f$ belong to the Besov space $B_{\be1}^1(\R^2)$
and let $N_1$ and $N_2$ be normal operators such that $N_1-N_2\in\bS_1$.
Then $f(N_1)-f(N_2)\in\bS_1$ and
$$
\|f(N_1)-f(N_2)\|_{\bS_1}\le\const\|f\|_{B_{\be1}^1(\R^2)}\|N_1-N_2\|_{\bS_1}.
$$
\end{thm}

\medskip

\section{\bf Operator H\"older functions and arbitrary moduli of continuity}

\medskip

For $\a\in(0,1)$, we consider the class $\L_\a(\R^2)$ of H\"older functions of order $\a$:
$$
\L_\a(\R^2)\df\left\{f:~\|f\|_{\L_\a(\R^2)}=
\sup_{z_1\ne z_2}\frac{|f(z_1)-f(z_2)|}{|z_1-z_2|^\a}<\be\right\}.
$$
The following result shows that in contrast with the class of Lipschitz functions, a H\"older function of order $\a\in(0,1)$ must be {\it operator H\"older of order $\a$}.

\begin{thm}
\label{oH}
There exists a positive number $c$ such that for every $\a\in(0,1)$ and every
$f\in\L_\a(\R^2)$,
$$
\|f(N_1)-f(N_2)\|\le c\,(1-\a)^{-1}\|f\|_{\L_\a(\R^2)}\|N_1-N_2\|^\a.
$$
for arbitrary normal operators $N_1$ and $N_2$.
\end{thm}

Consider now more general classes of functions. Let $\o$ be a modulus of continuity. We define the class $\L_\o(\R^2)$ by
$$
\L_\o(\R^2)\df\left\{f:~\|f\|_{\L_\o(\R^2)}=
\sup_{z_1\ne z_2}\frac{|f(z_1)-f(z_2)|}{\o(|z_1-z_2|)}<\be\right\}.
$$

As in the case of functions of one variable (see \cite{AP1}, \cite{AP2}), we define the function $\o_*$ by
$$
\o_*(x)\df x\int_x^\be\frac{\o(t)}{t^2}\,dt,\quad x>0.
$$

\begin{thm}
\label{Lo}
There exists a positive number $c$ such that for every
modulus of continuity $\o$ and every $f\in\L_\o(\R^2)$,
$$
\|f(N_1)-f(N_2)\|\le c\,\|f\|_{\L_\o(\R^2)}\,\o_*\big(\|N_1-N_2\|\big)
$$
for arbitrary normal operators $N_1$ and $N_2$.
\end{thm}

\begin{cor}
\label{sle}
Let $\o$ be a modulus of continuity such that
$$
\o_*(x)\le\const\o(x),\quad x>0,
$$
and let $f\in\L_\o(\R^2)$. Then
$$
\|f(N_1)-f(N_2)\|\le\const\|f\|_{\L_\o(\R^2)}\,\o\big(\|N_1-N_2\|\big)
$$
for arbitrary normal operators $N_1$ and $N_2$.
\end{cor}

\medskip

\section{\bf Perturbations of class $\bS_p$}

\medskip

In this section we study properties of $f(N_1)-f(N_2)$ in the case when $f\in\L_\a(\R^2)$, $0<\a<1$, and $N_1$ and $N_2$ are normal operators such that $N_1-N_2$
belongs to the Schatten--von Neumann class $\bS_p$. The following theorem generalizes Theorem 5.8 of \cite{AP3} to the case of normal operators.

\begin{thm}
\label{Sp}
Let $0<\a<1$ and $1<p<\be$. Then there exists a positive number $c$ such that
for every $f\in\L_\a(\R^2)$ and for arbitrary normal operators $N_1$ and $N_2$ with $N_1-N_2\in\bS_p$, the operator $f(N_1)-f(N_2)$ belongs to $\bS_{p/\a}$ and
the following inequality holds:
$$
\big\|f(N_1)-f(N_2)\big\|_{\bS_{p/\a}}\le c\,\|f\|_{\L_\a(\R^2)}\|N_1-N_2\|^\a_{\bS_p}.
$$
\end{thm}

For $p=1$ this is not true even for self-adjoint operators, see \cite{AP3}. Note that the construction of the counterexample in \cite{AP3} involves Hankel operators and is based on the criterion of membership of $\bS_p$ for Hankel operators, see \cite{Pe3}, Ch. 6.

The following weak version of Theorem \ref{Sp} holds:

\begin{thm}
\label{ws}
Let $0<\a<1$ and let $f\in\L_\a(\R^2)$. Suppose that $N_1$ and $N_2$ are normal operators such that $N_1-N_2\in\bS_1$. Then $f(N_1)-f(N_2)\in\bS_{\frac1\a,\be}$, i.e.,
$$
s_j\big(f(N_1)-f(N_2)\big)\le\const\|f\|_{\L_\a(\R^2)}(1+j)^{-\a},\quad j\ge0,
$$
\end{thm}

Here $s_j(T)$ is the $j$th singular value of a bounded operator $T$.

On the other hand, the conclusion of Theorem \ref{Sp} remains valid even for $p=1$ if we impose a slightly stronger assumption on $f$.

\begin{thm}
\label{s1a}
Let $0<\a<1$ and let $f$ belong to the Besov space $B_{\be1}^\a(\R^2)$. Suppose that $N_1$ and $N_2$ are normal operators such that $N_1-N_2\in\bS_1$. Then $f(N_1)-f(N_2)\in\bS_{1/\a}$ and
$$
\big\|f(N_1)-f(N_2)\big\|_{\bS_{1/\a}}\le\const\|f\|_{B_{\be1}^\a(\R^2)}
\|N_1-N_2\|_{\bS_1}^\a.
$$
\end{thm}

We conclude this section with the following improvement of Theorem \ref{Sp}.

\begin{thm}
\label{Spl}
Let $0<\a<1$ and $1<p<\be$. Then there exists a positive number $c$ such that
for every $f\in\L_\a(\R^2)$, every $l\in\Z_+$, and arbitrary normal operators $N_1$ and $N_2$ with bounded $N_1-N_2$, the following inequality holds:
$$
\sum_{j=0}^l\left(s_j\Big(\big|f(N_1)-f(N_2)\big|^{1/\a}\Big)\right)^p\le
c\,\|f\|_{\L_\a(\R^2)}^{p/\a}\sum_{j=0}^l\big(s_j(N_1-N_2)\big)^p.
$$
\end{thm}

\medskip

\section{\bf Commutators and quasicommutators}

We can generalize all the results listed above to the case of quasicommutators
$f(N_1)Q-Qf(N_2)$ and obtain estimates for such quasicommutators in terms of
$N_1Q-QN_2$ and $N_1^*Q-QN_2^*$. Here $Q$ is a bounded operator. In particular, if $Q$ is the identity operator, we arrive at the problem of estimating operator differences $f(N_1)-f(N_2)$ that has been considered above. On the other hand, if $N_1=N_2=N$, we arrive at the problem of estimating commutators $f(N)Q-Qf(N)$ in terms of $NQ-QN$ and $N^*Q-QN^*$.

\end{document}

%% file: CRinput.tex
\newcommand{\lb}{\linebreak}

\renewcommand{\a}{\alpha}

\newcommand{\z}{\zeta}

\newcommand{\s}{\sigma}

\renewcommand{\o}{\omega}

\renewcommand{\L}{\Lambda}

\newcommand{\F}{{\mathcal F}}

\newcommand{\X}{{\mathcal X}}

\newcommand{\C}{{\Bbb C}}

\newcommand{\R}{{\Bbb R}}
\newcommand{\Z}{{\Bbb Z}}

\newcommand{\bS}{{\boldsymbol S}}

\newcommand{\rf}[1]{(\ref{#1})}

\newcommand{\df}{\stackrel{\mathrm{def}}{=}}

\newcommand{\re}{\operatorname{Re}}

\newcommand{\const}{\operatorname{const}}

\newcommand{\eeq}{\end{equation}}
\newcommand{\beq}{\begin{equation}}
\newcommand{\bay}{\begin{eqnarray}}
\newcommand{\ba}{\begin{align*}}
\newcommand{\ea}{\end{align*}}
\newcommand{\ey}{\end{eqnarray}}
\newcommand{\bey}{\begin{eqnarray*}}
\newcommand{\eey}{\end{eqnarray*}}

\newcommand{\be}{\infty}

\newcommand{\im}{\operatorname{Im}}
\renewcommand{\re}{\operatorname{Re}}

\newtheorem{thm}{\hspace{\parindent}Theorem}[section]

\newtheorem{cor}[thm]{\hspace{\parindent}Corollary}